\documentclass[12pt,authoryear]{elsarticle}
\usepackage{graphicx}
\usepackage{mathrsfs}
\usepackage{amsmath}
\usepackage{amsfonts}
\usepackage[T2A]{fontenc}
\usepackage[english]{babel}
\usepackage{amssymb}%
\usepackage{natbib}
\usepackage{hyperref}
\newtheorem{theorem}{Theorem}
\newtheorem{acknowledgment}[theorem]{Acknowledgment}

\newtheorem{conclusion}[theorem]{Conclusion}

\newtheorem{definition}[theorem]{Definition}
\newtheorem{example}[theorem]{Example}

\newtheorem{remark}[theorem]{Remark}

\usepackage{setspace}
\author{ Martial Longla\\mlongla@olemiss.edu\\ Department of Mathematics, University of Mississippi, 38677 University, MS}\title{On mixtures of copulas and mixing coefficients}
\begin{document}

\begin{abstract} We show that if the density of the absolutely continuous part of a copula is bounded away from zero on a set of Lebesgue measure 1, then that copula generates \textquotedblleft lower $\psi$-mixing\textquotedblright\ stationary Markov chains. This conclusion implies $\phi$-mixing, $\rho$-mixing, $\beta$-mixing and \textquotedblleft interlaced $\rho$-mixing\textquotedblright . We also provide some new results on the mixing structure of Markov chains generated by mixtures of copulas.
\end{abstract}

\maketitle

\section{Introduction}
\subsection{Motivation and background}
The importance of mixing coefficients in the theory of central limit theorems for weakly dependent random sequences is an established fact. So, it is important to understand the mixing structure of various models. In the recent years, many researchers have been trying to provide sufficient conditions for mixing at various rates. These efforts have invited many people to investigate the properties of copulas, because they capture the dependence structure of stationary Markov chains. Mixtures of distributions are very popular in modeling. 

We present here a review of mixing coefficients and their convergence for mixtures of copulas. Convergence of mixing coefficients allow various central limit theorems for estimation of functions of the random variables or inference on model parameters. 

We have initially shown in \citet{LonPel2012}  a result on mixtures of copulas for absolute regularity. A second result on $\rho$-mixing for mixtures of copulas was provided in \citet{Longla2013}. Here, we extend the results to other mixing coefficients that are not less important.

We are providing here an improvement of one of our previous results. Namely, we have shown in Longla and Peligrad (2012), that Markov chains generated by a copula are $\phi$-mixing, when the density of its absolutely continuous part is bounded away from 0. This paper provides the proof that under this condition the Markov chains are lower $\psi$-mixing. Taking into account Theorem 1.2 and 1.3 of \citet{Bradley1997}, this conclusion implies \textquotedblleft interlaced $\rho$-mixing\textquotedblright. 

\subsection{Definitions}
Suppose $(\Omega,\mathscr{F},P)$ is a probability space, and $(X_{m}, m\in\mathbb{N})$ is a stationary Markov chain generated by a copula $C(x,y)$ on this space. Recall that a copula is a joint cumulative distribution function on $[0,1]^2$ with uniform marginals; and a copula-based Markov chain is nothing but a stationary Markov chain represented the copula of its consecutive states and its invariant distribution. Some examples of famous copulas are $\Pi(x,y)=xy$ - the independence copula, $W(x,y)=\max(x+y-1, 0)$- the Hoeffding lower bound and $M(x,y)=\min(x,y)$- the Hoeffding upper bound. $\Pi$ is the copula of independent random variables, $W$ induces on $[0,1]^2$ a singular probability measure with support $S_1=\{ (x,y)\in[0,1]^2 : x+y=1\}$ and $M$ induces a singular probability measure with support $S_2=\{ (x,y)\in[0,1]^2 : x=y\}$. The copula $C(x,y)$ induces on the unit square a probability measure that we will denote $P_C$. This probability measure acts on sets $A=(0, x]\times(0, y]$ by $P_C (A)=C(x,y)$. For a set of copulas $C_1 (x,y), \cdots, C_k (x,y)$ and some strictly positive numbers $a_1 ,\cdots, a_k$ such that $\displaystyle \sum_{i=1}^{k}a_i = 1$, $\displaystyle C(x,y)=\sum_{i=1}^{k}a_i C_i (x,y) $ is a copula and the measure it induces on $[0,1]^2$ is $\displaystyle P_C =\sum_{i=1}^{k}a_i P_{C_i}$. Let also $P_{CF}$ denote the probability measure induced by $C(x,y)$ and the univariate cumulative distribution function $F$ on our probability space. Let $f_{,i}(x,y)$ denote the derivative with respect to the $i^{th}$ variable of the function $f$ at the point $(x,y)$. The fold product of the copulas $C_1 (x,y)$ and $C_2 (x,y)$ is the copula defined by 

$$ C(x,y)=C_{1}*C_{2}(x,y)=\int_{0}^{1}C_{1,2}(x,t)C_{2,1}(t,y)dt.$$ 

\noindent This operation is associative and distributive over convex combinations of copulas. The widely used notation of the powers of $C(x,y)$ is $C^{m}(x,y)=C^{m-1}*C(x,y)$ with $C^{1}(x,y)=C(x,y)$ is called $m^{th}$ fold product of $C(x,y)$. This copula represents that of the random vector $(X_{0}, X_{m})$, where $(X_{0}, X_{1}, \cdots, X_{m})$ is a stationary Markov chain with copula $C(x,y)$ and the uniform marginal.

Note that for any vector random variable $(U,V)$ with joint distribution $H(u,v)$ and respective marginals $F_U$ and $F_V$, $H(u,v)=C(F_U (u), F_V (v))$ for some copula $C(x,y)$, where $F_X$ represents the distribution of the random variable $X$. Sklar's theorem ensures uniqueness of this representation for continuous random variables. Let $c(x,y)$ be the density of the absolutely continuous part of the copula $C(x,y)$ and $c_{m}(x,y)$ the density of the absolutely continuous part of the $m^{th}$ fold product of $C(x,y)$. For definitions, check \citet{LonPel2012} or \citet{Longla2013}. 

Using the Lebesgue measure as invariant distribution after rescaling the effect of the marginal distribution, taking into account the notations from \citet{LonPel2012} the mixing coefficients of interest for an absolutely continuous copula and an absolutely continuous invariant distribution for the states of the stationary Markov chain that they generate are defined below.

\begin{definition} Let $\mathcal{B}$ denote the Borel $\sigma$-algebra on $[0,1]$. Let $\lambda$ be the Lebesgue measure on $[0,1]$, $B^{c}$ be the complement of the set $B$ in $[0,1]$. $\mathbb{E}_{\mu}(f)$ is the expected value of the random variable $f(X)$ under the probability measure $\mu$ and $\displaystyle ||f||_2=\Big(\int_{0}^{1}f^2 (x)dx\Big)^{1/2}. $ \\ 
\noindent $\displaystyle\rho_n=\sup_{f,g}\Big\{\int_{0}^{1}\int_{0}^{1}c_{n}(x,y)f(x)g(y)dxdy, \int_{0}^{1}f(x)dx=0,\int_{0}^{1}g(x)dx=0, ||f||_2 = ||g||_2 = 1 \Big\} $ is the maximal coefficient of correlation. The Markov chain is $\rho$-mixing, if $\rho_{n}$ converges to 0.\\ $\displaystyle\phi_{n}=\sup_{B\subset\mathcal{B}\cap[0,1]}ess\sup_{x\in[0,1]}|\int_{B}(c_{n}(x,y)-1)dy|$ is the uniform mixing coefficient. The Markov chain is uniformly mixing, if $\phi_{n}$ converges to 0. \\ $\displaystyle\beta_n=\int_{0}^{1}\sup_{B\in\mathcal{B}\cap[0,1]}|\int_{B}(c_{n}(x,y)-1)dy|dx$ is the absolute regularity mixing coefficient. The Markov chain is absolutely regular, if $\beta_{n}$ converges to 0. If the convergence of $\beta_{n}$ is exponential, then we say that the Markov chain is geometrically ergodic.\\ $\displaystyle\psi^{\prime}_{n}=\inf_{A,B\in\mathcal{B}, \lambda(A)\lambda(B)>0}\frac{\int_{A}\int_{B}c_{n}(x,y)dxdy}{\lambda(A)\lambda(B)}$ is the \textquotedblleft lower $\psi$-mixing\textquotedblright. We say that the Markov chain is \textquotedblleft lower $\psi$-mixing\textquotedblright\ or $\psi^{\prime}$-mixing, if $\psi^{\prime}_{n}$ converges to 1.\\ $\displaystyle\psi_{n}=\sup_{A,B\in\mathcal{B}, \lambda(A)\lambda(B)>0}\frac{|\int_{A}\int_{B}(c_{n}(x,y)-1)dxdy|}{\lambda(A)\lambda(B)}$ is the $\psi$-mixing coefficient. The Markov chain is $\psi$-mixing if $\psi_{n}$ converges to 0.
\end{definition} 

These coefficients have a more complex formulation for a general random sequence. For the sake of clarity, we shall provide the general definitions of these coefficients for stationary Markov chains with marginal $F$. Let $\mathcal{R}=\sigma(X_i)$ be the $\sigma$-algebra generated by $X_i$ and $\mathbb{L}_{2}(\mathcal{\sigma}(X_i ))$ be the space of random variables that are $\mathcal{\sigma}(X_i ) $ measurable and square integrable. Let $Corr(X,Y)$ be the correlation coefficient of the random variables $X$ and $Y$. The derivation of the above equivalent forms is provided in \citet{LonPel2012}.

\begin{definition} Under the above assumptions and notations, let $\mathcal{R}^2=\sigma(X_0, X_n)$,  where $\sigma(X_{m}, m\in S)$ is the $\sigma$-algebra generated by the random variables indexed by $S$. For every positive integer $n$, let $\mu_n$ be the measure induced by the distribution of $(X_0, X_n)$. Let $\mu$ be the measure induced by the distribution of $X_{0}$. The coefficients of interest are defined as follows:

\[\displaystyle
\rho_n=\rho (\sigma(X_0 ), \sigma(X_n )):=\sup\{Corr(f,g),\text{ }f\in\mathbb{L}%
_{2}(\mathcal{\sigma}(X_0)),\text{ }g\in\mathbb{L}_{2}(\mathcal{\sigma}(X_n))\},
\]

\[\displaystyle
\phi_n=\phi (\sigma(X_0 ), \sigma(X_n )):=\sup_{A,B\in\mathcal{R}%
, \mu(A)>0}|P(X_n \in B | X_0 \in A)-\mu(B)|,
\]

\[\displaystyle
\beta_n=\beta (\sigma(X_0 ), \sigma(X_n )):=\sup_{D\in \mathcal{R}^2}|\mu_n (D)-(\mu\times \mu)(D)|,
\]

\[\displaystyle
\psi^{'}_n=\beta (\sigma(X_0 ), \sigma(X_n )):=\inf_{A,B\in\mathcal{R}, \mu(A)\mu(B)>0}\frac{\mu_n (A\times B)}{\mu(A)\mu(B)},
\]

\[\displaystyle
\psi_n=\psi (\sigma(X_0 ), \sigma(X_n )):= \sup_{ A,B\in\mathcal{R}, \mu(A)\mu(B)>0}\frac{|\mu_n (A\times B)-\mu(A)\mu(B)|}{\mu(A)\mu(B)},
\]

\noindent \leftline{$\displaystyle\rho_{n}^{*}=\sup\Big\{|Corr(f,g)|, f\in \mathbb{L}_{2}(\sigma(X_{i}, i\in S)), g\in \mathbb{L}_{2}(\sigma(X_{i}, i\in T)), S\subset \mathbb{Z}, T\subset\mathbb{Z}$,} \centerline{ $\displaystyle \min_{s\in S, t\in T}|s-t|\geq n\Big\}.$}

\end{definition}

$\rho^{*}_n$ is the \textquotedblleft interlaced $\rho$-mixing\textquotedblright\ coefficient. The Markov chain is \textquotedblleft interlaced $\rho$-mixing\textquotedblright, if $\rho^{*}$ converges to 0.

\subsection{Setup and Results}
The main purpose of this paper is to prove three main results, stated as Theorems \ref{psi}, \ref{mixture} and \ref{mult} below. The background and the proof of those results are given in Section 2. Section 3 gives examples of well-known copulas to which these results apply. 

\begin{theorem}[Lower $\psi$-mixing]\label{psi}
Let $(X_{i}, i\in \mathbb{N})$ be a stationary Markov chain generated by a copula $C(x,y)$ and an absolutely continuous marginal distribution $F$ on $\mathbb{R}$. Suppose that for some positive integer $m$, the copula of the random vector $(X_0, X_m )$ has an absolutely continuous component with density function $c_{m}(x,y)$ bounded away from zero on a set of Lebesgue measure 1. Then the Markov chain is \textquotedblleft lower $\psi$-mixing\textquotedblright\ (with $1-\psi^{'}_n \to 0$ at least exponentially fast as $n\to\infty$ ) and hence also $\phi$-mixing, $\rho$-mixing, $\beta$-mixing and \textquotedblleft interlaced $\rho$ -mixing\textquotedblright\ (all with mixing rates at least exponentially fast).
\end{theorem}

\begin{theorem}[Mixtures]\label{mixture}
Let $A=(X_{i}, i\in \mathbb{N})$ be a Stationary Markov chain generated by a convex combination of copulas $C_{i}(x,y), 1\leq i\leq k$ and the uniform distribution on $[0,1]$. Suppose $m$ is a positive integer. The following statements hold:
\begin{enumerate}
\item If one of the copulas in the combination generates stationary Markov chains satisfying $\rho_m <1$, then the Markov chain $A$ satisfies $\rho_m<1$ and is $\rho$-mixing with at least exponential mixing rate.
\item If one of the copulas in the combination generates stationary Markov chains satisfying $\psi^{\prime}_m>0$, then the Markov chain $A$ satisfies $\psi^{\prime}_m>0$ and is $\psi^{\prime}$-mixing with at least exponential mixing rate.
\item If one of the copulas in the combination generates stationary $\psi$-mixing Markov chains, then the Markov chain $A$ is $\psi^{\prime}$-mixing with at least exponential mixing rate.
\item If one of the copulas in the combination generates ergodic and aperiodic stationary Markov chains satisfying $\phi_{m}<1$, then the Markov chain $A$ satisfies $\phi_m$ and is $\phi$-mixing with at least exponential mixing rate.
\item If one of the copulas in the combination generates ergodic and aperiodic stationary Markov chains with $\beta_{m}<1$, then the Markov chain $A$ satisfies $\beta_m<1$ and is $\beta$-mixing.

\end{enumerate} 
\end{theorem}
Statement 3 is an immediate corollary of Statement 2, since $\psi$-mixing implies $\psi^{'}$-mixing. Statement 3 emphasizes the fact that even if one of the components of the copula generates a stationary $\psi$-mixing Markov chain, the Markov chain $A$ itself is not necessarily $\psi$-mixing - one can only conclude that it is $\psi^{'}$-mixing. Statements 1, 2, 4, and 5 are special cases of the corresponding parts of Theorem \ref{mult} below. 
Statements 1, 2, 4 and 5 of Theorem \ref{mixture} will be strengthened in Theorem \ref{mult} below via the following setup. 

Suppose $k$ is a positive integer, for each $i\in\{1,\cdots, k\}$, $C_{i}(x,y)$ is a copula, $a_1, \cdots, a_k$ are strictly positive integers such that $a_1+\cdots+a_k =1$, and $C(x,y)$ is the copula defined by 

$$\displaystyle C(x, y)=\sum_{i=1}^{k}a_i C_i (x,y).$$

Suppose $m$ is a positive integer. For each $m$-tuple $(i)=(i_1, \cdots, i_m)$ of elements of $\{1, \cdots, k\}$, let $(X_0^{(i)},\cdots, X_m^{(i)})$ be a Markov random vector such that for each $j\in \{1,\cdots,m\}$, the joint distribution of the random vector $(X_{j-1}^{(i)}, X_{j}^{(i)})$ is $D^{(i)}_j (x,y)=C_{i_j}(x,y)$. 
Using the Chapman-Kolmogorov equations and properties of copulas we find the cumulative distribution function of $(X^{(i)}_{0},\cdots, X^{(i)}_{m} )$ as 

$$\displaystyle F^{(i)}(t_0, \cdots, t_m)=\nu^{(i)}((0,t_0]\times\cdots\times(0, t_m] )= \int_{0}^{t_m}\cdots\int_{0}^{t_0}dx_0 C_{i_1}(x_0, dx_1 )\cdots C_{i_m}(x_{m-1}, dx_m ),$$

\noindent where $\nu^{(i)}$ is the probability measure induced by $F^{(i)}$ on $\displaystyle [0,1]^{m+1}$ and for $i\in\{0,\cdots, m\}$, $t_i\in [0,1]$. On the other hand, the cumulative distribution function of $(X_{0},\cdots, X_{m} )$ is 

$$\displaystyle F(t_0, \cdots, t_m)= \int_{0}^{t_m}\cdots\int_{0}^{t_0}dx_0 C(x_0, dx_1 )\cdots C(x_{m-1}, dx_m )=$$ $$=\sum_{(i)}a^{(i)}\int_{0}^{t_m}\cdots\int_{0}^{t_0}dx_0 C_{i_1}(x_0, dx_1 )\cdots C_{i_m}(x_{m-1}, dx_m ),$$ 

\noindent where $\displaystyle a^{(i)}=\Pi_{j=1}^{m}a_{i_j}$. The last equality uses linearity of the integral and the special form of $C(x,y)$. So, if $\nu$ is the probability measure induced by $F$ on $[0,1]^{m+1}$, we obtain for every Borel measurable $D$ of $[0,1]^{m+1}$, 
\begin{equation}\displaystyle \nu(D)=\sum_{(i)}a^{(i)}\nu^{(i)} (D).\label{nu} \end{equation} 

Now, consider $\mu_m$ - joint distribution of $(X_0, X_m)$. For each $m$-tuple $(i)$ of elements of $\{1,\dots,k\}$, let $\mu_m^{(i)}$ denote the distribution of the random vector $(X_0^{(i)}, X_m^{(i)})$. Then, as a special case of equation (\ref{nu}), one has that for every set $Q\in \mathcal{R}^{2}$, 

\begin{equation}\displaystyle \mu_{m}(Q)=\sum_{(i)}a^{(i)}\mu_m^{(i)} (Q).\label{mu}\end{equation} 

Formula (\ref{nu}) is key to the proof of the following.

\begin{theorem}[Extensions]\label{mult} Consider the convex combination copula $C(x,y)$ in the above setup. Suppose $m$ is a positive integer. Consider the above notations, including the Markov random vector $(X_0^{(i)},\cdots, X_m^{(i)})$. The following statements hold:

\begin{enumerate}

\item If for some $m$-tuple $(i_0 )$ (of elements of $\{1,\cdots, k\}$), one has that $\rho(\sigma(X_{0}^{(i_0 )}), \sigma(X_{m}^{(i_0 )}))<1$, then the Markov chain $(X_{i}, i\in\mathbb{N})$ generated by the convex combination $C(x,y)$ and the uniform distribution on $[0,1]$ satisfies $\rho_m<1$ and is exponential $\rho$-mixing.

\item If for some $m$-tuple $(i_0 )$ (of elements of $\{1,\cdots, k\}$), one has that $\psi^{\prime}(\sigma(X_{0}^{(i_0 )}), \sigma(X_{m}^{(i_0 )}))>0$, then the Markov chain $(X_j, j\in\mathbb{N})$ generated by the convex combination $C(x,y)$ satisfies $\psi^{'}_m>0$ and is exponential $\psi^{\prime}$-mixing.

\item Suppose that one has that at least one of the copulas in the combination generates a stationary, ergodic and aperiodic Markov chain. If for some $m$-tuple $(i_0 )$ (of elements of $\{1,\cdots, k\}$), one has that $\phi (\sigma(X_{0}^{(i_0 )}), \sigma(X_{m}^{(i_0 )}))<1$, then the Markov chain $(X_j, j\in\mathbb{N})$ generated by the convex combination $C(x,y)$ and the uniform distribution on $[0,1]$ satisfies $\phi_m<1$ and is exponential $\phi$-mixing.

\item Suppose that one has that at least one of the copulas in the combination generates a stationary, ergodic and aperiodic Markov chain. If for some $m$-tuple $(i_0 )$ (of elements of $\{1,\cdots, k\}$), one has that $\beta(\sigma(X_{0}^{(i_0 )}), \sigma(X_{m}^{(i_0 )}))<1$, then the Markov chain $(X_j, j\in\mathbb{N})$ generated by the convex combination and the uniform distribution on $[0,1]$ satisfies $\beta_m<1$ and is absolutely regular.

\end{enumerate}

\end{theorem}

\section{Proofs of the Theorems}

\subsection{Theorem \ref{psi}}
The condition of Theorem \ref{psi} was first investigated by \citet{Beare2010}, who showed that it implies $\rho$-mixing for $m=1$. Then, \citet{LonPel2012} have shown that the condition implies $\phi$-mixing for $m=1$. During a discussion at the conference on Stochastic Processes and Application (SPA 2013), Bradley suggested that this condition could imply \textquotedblleft lower $\psi$-mixing\textquotedblright. The following well-known result is provided in Remark 1.1 of \citet{Bradley1997} and is used in the proofs.

\begin{remark}\label{remark}
A stationary Markov chain is exponentially $\psi^{\prime}$-mixing, if for some integer $m$, $\psi^{\prime}_{m}>0$. 
\end{remark}

\citet{Bradley1997} also showed, based on this result, that if for some $m$, $\psi^{\prime}_{m}>0$, then the Markov chain is $\rho^{*}$-mixing with exponential decay rate. We will use this result in the context of copula-based Markov chains to establish the proof of Theorem \ref{psi}.

Let $(X_n, n\in\mathbb{N})$ be a Markov chain generated by the uniform distribution on $[0,1]$ and a copula $C(x,y)$ with density of the absolutely continuous part of $C^{m}(x,y)$ greater than $c>0$ on a set of Lebesgue measure 1 for some integer $m\ge 1$. We shall first show that if $c_m (x,y)\ge c \quad a.e.$ on $[0,1]^2$, then 
\begin{equation} \label{bound}
P(X_0 \in A , X_m \in B)\ge c\lambda(A)\lambda(B), \quad \mbox{for all} \quad A, B\in\mathcal{B}.
\end{equation} 
Recall that for any copula $C(x,y)$, there exists a unique representation $$C(x,y)=AC(x,y)+SC(x,y),$$ where $AC(x,y)$ is the absolutely continuous part of $C(x,y)$. $AC$ induces on $[0,1]^{2}$ a measure $P_{c}$ defined on Borel sets by $$\displaystyle P_{c}(A\times B)=\int_{A}\int_{B}c_m (x,y)dydx.$$ $SC(x,y)$ is the singular part of the copula. It induces a singular measure on $[0,1]^{2}$. If we keep the notation $SC$ for this singular measure, then $\displaystyle P(X_0 \in A , X_m \in B)=P_{c}(A\times B)+SC(A\times B).$ Thus, 

\begin{equation}\displaystyle P((X_0, X_m)\in A\times B)\ge P_{c}(A\times B)=\int_{A}\int_{B}c_m (x,y)dydx.\label{bound1}\end{equation}

Taking into account the fact that $c_m (x,y)\ge c \quad a.e. $, we obtain $$P(X_0 \in A, X_m \in B)\ge c\lambda(A)\lambda(B) \mbox{ for all } A, B\in \mathcal{B}.$$ Therefore, inequality (\ref{bound}) holds.

A similar formula will hold for any Markov chain generated by $C(x,y)$ and any absolutely continuous marginal distribution $F(x)$ of $X_0 $. This is due to the fact that, if $f(x)$ is the density of $F(x)$, then $$c_m(F(x), F(y) )\geq c $$ a.e. by the assumption. Thus, the joint density of $(X_0,X_m)$ - $h_m(x,y)$ will exist and satisfy almost everywhere $$h_m(x,y)=c_m(F(x),F(y))f(x)f(y)\geq cf(x)f(y),$$ where $f$ is the density of $X_0$. One will just need to replace $c_m(x,y)$ by $h_m(x,y)$ in inequality (\ref{bound1}). This inequality leads to $\psi^{\prime}_{1}\ge c>0$ for the considered Markov chain. 

Therefore, by Remark \ref{remark}, the Markov chains generated by $C(x,y)$ and any absolutely continuous univariate marginal distribution are $\psi^{\prime}$-mixing and $\rho^{*}$-mixing. It is also well known that $\psi^{\prime}$-mixing implies $\phi$-mixing, which implies $\rho$-mixing and exponential $\beta$-mixing, confirming the results of previous studies.

\subsection{Theorem \ref{mult}} 
{\quad}
Let $\displaystyle C(x,y)=\sum_{k=1}^{k}a_{i}C_{i}(x,y)$, $a_{i}> 0$, $\displaystyle \sum_{i=1}^{k}a_{i}=1$ and assume in each of the cases that $C_{1}(x,y)$ is the copula that generates stationary Markov chains with the given mixing property.

The proof of Statement 1 of Theorem \ref{mixture} is similar to the one provided for Theorem 5 in \citet{LonPel2012}. It is solely based on the fact that we need only to show that $\rho_{m}<1$. The fact that the coefficients of the convex combination add up to 1 and $\rho_{n}\leq 1$ helps in the conclusion. In fact, $\rho_{m}<1$ for some $m$ iff $\rho_{m}$ converges to 0. This fact can be easily verified by looking at the sequence $\rho_{im}\leq \rho_{m}^i$ as $i\to\infty$ for a fixed value of $m$, and using $\rho_{n+1}\leq \rho_{n}$.
Let $f$ and $g$ defined on $[0,1]$ and satisfy $\displaystyle \mathbb{E}_\lambda (f(X_m))=0, $ $\mathbb{E}_\lambda (g(X_0))=0$, and $\displaystyle \displaystyle \mathbb{E}_\lambda (f^2 (X_m))=\mathbb{E}_\lambda (g^2 (X_0))=1$. Under these conditions, $\displaystyle Corr(g(X_0), f(X_m))=\mathbb{E}_{\mu_m} (g(X_0)f(X_m))$. The fold product of copulas yields $$\displaystyle \rho_m = \sup\{\int_{0}^{1}\int_{0}^{1}f(x)g(y)C^{m}(dx,dy), \quad \int_{0}^{1}f(x)dx=\int_{0}^{1}g(y)dy=0 , ||f||_2 = ||g||_2 =1 \}.$$ Using formula (\ref{nu}) and properties of expectations, we obtain 

$$\displaystyle Corr (g(X_0), f(X_m))=\sum_{(i)}a^{(i)}\mathbb{E}_{\mu_m^{(i)}}(g(X^{(i)}_0) f(X^{(i)}_m)),$$ 

where the sum is taken over all possible $(i)$. Note that under these assumptions, $\displaystyle \mathbb{E}_\lambda (f(X^{(i)}_m))=0$, $\mathbb{E}_\lambda (g(X^{(i)}_0))=0$, and $\displaystyle \displaystyle \mathbb{E}_\lambda (f^2 (X^{(i)}_m))=\mathbb{E}_\lambda (g^2 (X^{(i)}_0))=1$. So, 

$$\displaystyle Corr (g(X_0), f(X_m))=\sum_{(i)}a^{(i)}corr(g(X^{(i)}_0), f(X^{(i)}_m)).$$ 

Therefore, $\displaystyle \rho (\sigma(X_0 ), \sigma(X_m ))\leq \sum_{(i)}a^{(i)}\rho(\sigma (X^{(i)}_0), \sigma (X^{(i)}_m))$. 

\noindent Provided that $\displaystyle \rho(\sigma (X^{(i_{0} )}_0), \sigma (X^{(i_{0} )}_m))<1$ for this value of $m$ and some ($i_0$), we obtain 

$$\displaystyle \rho (\sigma(X_0 ), \sigma(X_m ))\leq a^{(i_{0})}\rho(\sigma (X^{(i_{0})}_0), \sigma (X^{(i_{0})}_m)) -a^{(i_{0})}+ a^{(i_{0})}+\sum_{(i)\neq (i_{0})}a^{(i)}\rho(\sigma (X^{(i)}_0), \sigma (X^{(i)}_m)).$$

Given that $\displaystyle \rho(\sigma (X^{(i)}_0), \sigma (X^{(i)}_m))\leq 1$, it follows that 

$$\displaystyle \rho (\sigma(X_0 ), \sigma(X_m ))\leq a^{(i_{0})}(\rho(\sigma (X^{(i_{0})}_0), \sigma (X^{(i_{0})}_m)) -1)+\sum_{(i)}a^{(i)}=1-(1-\rho(\sigma (X^{(i_{0})}_0), \sigma (X^{(i_{0})}_m)))a^{(i_{0})}.$$

This inequality takes into account the fact that $\displaystyle \sum_{(i)}a^{(i)}=(a_1 + \cdots+a_k)^{m}=1$ and $a_i> 0$. Therefore, $\rho_m<1$, indicating that the Markov chain generated by the convex combination is $\rho$-mixing.

The proof of Statement 2 of Theorem \ref{mult} follows from Remark \ref{remark} and the fact that the given convex combination generates stationary Markov chains for which equation $(\ref{mu})$ holds.

\noindent Note that if $\displaystyle \psi^{\prime}_{*m}=\psi^{'}(\sigma(X_0^{(i_0)}), \sigma(X_{m}^{(i_0)}))$ and $\displaystyle C(x,y)=\sum_{i=1}^{k}a_{i}C_{i}(x,y)$ with $0< a_i\leq 1$ for $i\ge 1$, then $\displaystyle\psi^{\prime}_{*m}>0$ by our assumptions. Now, we bound the $m^{th}$ mixing coefficient for the mixture using equation (\ref{mu}). For any $A, B \in \mathcal{B}$ such that $\lambda(A)\lambda(B)>0$,

$$\frac{P(X_0\in A, X_m\in B)}{P(X_0\in A)\cdot P(X_m\in B)}=\frac{\mu_m(A\times B)}{\lambda(A)\cdot \lambda(B)}=\sum_{(i)}a^{(i)}\frac{\mu_m^{(i)}(A\times B)}{\lambda(A)\cdot\lambda(B)}\geq a^{(i_0)}\frac{\mu_m^{(i_0)}(A\times B)}{\lambda(A)\cdot\lambda(B)}.$$ 

It follows that 

$$\displaystyle\psi^{\prime}_{m}= \inf_{A,B\in\mathcal{B}, \lambda(A)\lambda(B)>0}\frac{\mu_m( A\times B)}{\lambda(A)\cdot\lambda( B)}\ge a^{(i_0)}\inf_{A,B\in\mathcal{B}, \lambda(A)\lambda(B)>0}\frac{\mu_m^{(i_0)}(A\times B)}{\lambda(A)\cdot \lambda(B)} =a^{(i_0)}\psi^{\prime}_{*m}>0.$$

The proofs of Statement 3 and 4 of Theorem \ref{mult} use equation (\ref{mu}) and the fact that if at least one of the copulas in the combination generates a stationary Markov chain that is ergodic and aperiodic, then the Markov chain itself is ergodic and aperiodic. The later fact was proved in Lemma 3 of \citet{LonPel2012}.

Under the assumptions of Statement 3 of Theorem \ref{mult}, let $C_{1}(x,y)$ generate a stationary ergodic and aperiodic Markov chain. It follows that the stationary Markov chain generated by any mixture of copulas that contains $C_{1}(x,y)$ is aperiodic and ergodic. 

Therefore, the Markov chain generated by convex combinations containing $C_1 (x,y)$ is $\phi$-mixing if $\phi_{a}<1$ for some positive integer $a$, by Theorem 21.22 of \citet{Bradley2007}. 

Now, let $\phi_{*m}= \phi_m(\sigma(X_{0}^{(i_0)}), \sigma(X_m^{(i_0)}))$. For this value of $m$ and any $A,B\in\mathcal{B}$ such that $\lambda(A)>0,$ by equation (\ref{mu}) and the fact that $\displaystyle \sum_{(i)}a^{(i)}=1,$ 

$$|P(X_{m}\in B| X_0\in A)-P(X_m\in B)|=\frac{|\mu_m (A\times B)-\lambda(A)\cdot\lambda(B)|}{\lambda(A)} =$$ 
$$=\Big|\sum_{(i)}a^{(i)}\frac{\mu_{m}^{(i)}(A\times B)-\lambda(A)\cdot\lambda(B)}{\lambda(A)}\Big|\leq \sum_{(i)}a^{(i)}\frac{|\mu_m^{(i)}(A\times B)-\lambda(A)\cdot\lambda(B)|}{\lambda(A)}.$$
Therefore, using the fact that $\phi$-mixing coefficients are less than or equal to 1, it follows that
$$\phi(\sigma(X_0), \sigma(X_m))\leq \sum_{(i)}a^{(i)}\phi(\sigma(X_0^{(i)}), \sigma(X_m^{(i)}))\leq a^{(i_0)}\phi(\sigma(X_0^{(i_0)}), \sigma(X_m^{(i_0)}))+\sum_{(i)\neq (i_0)}a^{(i)}.$$ 
So, $\displaystyle\phi(\sigma(X_0), \sigma(X_m))\leq a^{(i_0)}( \phi_{*m}-1)+1 <1$. 

Therefore, the Markov chain generated by the convex combination $C(x,y)$ containing $C_{1}(x,y)$ and any absolutely continuous univariate marginal distribution are $\phi$-mixing under the assumptions of Statement 3 of Theorem \ref{mult}.

Taking into account the above comments, Theorem 21.5 and Corollary 21.7 in \citet{Bradley2007}, for the proof of Statement 4 of Theorem \ref{mult}, it is enough to show that $\beta(\sigma(X_0), \sigma(X_m))<1$.
Let $D\in \mathcal{B}^2$. By equation (\ref{mu}) and the fact that $\displaystyle \sum_{(i)}a^{(i)}=1$,

$$|\mu_m(D)-(\lambda\times\lambda)(D)|=\Big|\sum_{(i)}a^{(i)}(\mu_m^{(i)}(D)-(\lambda\times\lambda)(D))\Big|\leq \sum_{(i)}a^{(i)}|\mu_m^{(i)}(D)-(\lambda\times\lambda)(D)|.$$
Therefore, using the fact that $\beta$-mixing coefficients are less than or equal to 1, 

$$\beta(\sigma(X_0), \sigma(X_m))\leq \sum_{(i)}a^{(i)}\beta(\sigma(X^{(i)}_0), \sigma(X^{(i)}_m))\leq a^{(i_0)}\beta(\sigma(X^{(i_0)}_0), \sigma(X^{(i_0)}_m))+ \sum_{(i)\neq (i_0)}(a^{(i)}\cdot 1).$$
It follows that $\displaystyle \beta(\sigma(X_0), \sigma(X_m))\leq a^{(i_0)}(\beta(\sigma(X^{(i_0)}_0), \sigma(X^{(i_0)}_m))-1)+1<1$. The result follows.

\section{Examples and remarks}
\begin{example}
Any convex combination of copulas that contains the independence copula has the density of its absolutely continuous part bounded away from 0, and Theorem \ref{psi} applies to it.
\end{example}
 
In fact, for any convex combination of copulas containing the independence copula, the density of the absolutely continuous part is greater than the coefficient of the independence copula in the combination.
 
 An example of copula family of the above kind is the Frechet family of copulas defined as follows: $$\displaystyle C_{a,b}(x,y)=aW(x,y)+bM(x,y)+(1-a-b)\Pi(x,y), a\geq 0, b\geq 0, a+b\leq 1.$$ 
A subset of the Frechet family of copulas is the Mardia family of copulas, defined by 
$$\displaystyle C_{\theta}(x,y)=\frac{\theta^2 (1+\theta)}{2}M(x,y)+(1-\theta^2)\Pi(x,y)+\frac{\theta^2 (1-\theta)}{2}W(x,y), | \theta | \leq 1,$$  if we take $a+b=\theta^2$. Copulas from the Frechet family satisfy Theorem \ref{psi} for $a+b\neq 1$.

 These examples where shown to generate $\phi$-mixing in \citet{LonPel2012}. 

To emphasize the importance of Statement 3 of Theorem \ref{mixture}, consider the stationary Markov chain $(X_{i}, i\in\mathbb{N})$ generated by a copula from the above families. The copula of $(X_{0}, X_{n})$ is a member of the family. Denote this copula $$\displaystyle C^{n}(x,y):=a_n W(x,y)+b_n M(x,y)+(1-a_n -b_n)\Pi(x,y).$$ 

Let $F$ be the marginal cumulative distribution function of $X_0$ and $\mu$ the measure induced by $F$. Assume that $F$ is absolutely continuous and strictly increasing. We shall show that this Markov chain is not $\psi$-mixing. For any $A,B\in\mathcal{R},$ we have 
$$\displaystyle P(X_0 \in A , X_n\in B)=P_{C^{n}F}(A\times B)=a_n P_{WF}(A\times B)+b_n P_{MF}(A\times B)+(1-a_n -b_n)\mu(A)\mu(B).$$ 

Thus, $$\displaystyle \frac{\mu_n(A \times B)-\mu(A)\mu(B)}{\mu(A)\mu(B)}=a_n \frac{P_{WF}(A\times B)-\mu(A)\mu(B)}{\mu(A)\mu(B)}+b_n \frac{P_{MF}(A\times B)-\mu(A)\mu(B)}{\mu(A)\mu(B)}.$$ 

Take sets $A,B\in\mathcal{R}$ as follows: $A=B=F^{-1}(1/2-\frac{\varepsilon}{2},1/2+\frac{\varepsilon}{2})$. Clearly, $\mu(A)=\mu(B)= \varepsilon$ and $P_{WF}(A\times B)=P_W \Big((1/2-\frac{\varepsilon}{2},1/2+\frac{\varepsilon}{2})^2\Big) = \varepsilon$, and  $P_{MF}(A\times B)=P_M \Big((1/2-\frac{\varepsilon}{2},1/2+\frac{\varepsilon}{2})^2\Big)= \varepsilon$. 
Combining this with the definition of $\psi_n$ leads to
$$ \psi_n \geq \frac{|\mu_n(A\times B)-\mu(A)\mu(B)|}{\mu(A)\mu(B)}= \frac{(1-\varepsilon)(a_n+b_n)}{\varepsilon}.$$ 

Therefore, $\displaystyle\psi_n\geq\frac{(1-\varepsilon)(a+b)^n}{\varepsilon}$, using \citet{Longla2014}. It follows that, if the distribution of $X_0$ is absolutely continuous and strictly increasing, then taking the supremum over all such $A,B\in\mathcal{R}$ for a fixed integer $n\ge 1$ leads to $\psi_n =\infty$. 
This clarifies the following remark. 

\begin{remark}
Copulas from the Mardia and Frechet families don't generate $\psi$-mixing continuous space strictly stationary Markov chains for any $a,b$. So, \textquotedblleft lower $\psi$-mixing\textquotedblright\ is the best result we can obtain for a Markov chain generated by a copula from these families. Notice that the marginal distribution plays no role in this result. This fact also shows that Statement 3 of Theorem \ref{mixture} cannot be strengthened to a conclusion on $\psi$-mixing for the convex combinations of copulas in general. Also, $\psi$-mixing does not follow in general from the assumptions of Theorem \ref{psi}. We can also mention that the results of Theorems \ref{mixture} and \ref{mult} do not require the marginal distribution to be uniform. These results will still be valid for any marginal distribution.
\end{remark}
\begin{example}
The Marshall-Olkin family of copulas defined by $\displaystyle C_{a,b}(x,y)=\min(xy^{1-a},yx^{1-b})$ with $0\leq a, b\leq 1$ satisfies Theorem \ref{psi} for values of $a, b$ such that $a\neq 1$ and $b\neq 1$. 
\end{example}
In fact, for copulas from this family, we have $c_{a,b}(x,y)\ge \min(1-a,1-b)$, except for $y^a=x^{b}$, which is a set of Lebesgue measure 0. \citet{LonPel2012} have shown that this family of copulas generates $\phi$-mixing, which is a weaker result than what we have shown here because $\psi^{'}$-mixing implies $\phi$-mixing.

\begin{conclusion}

It is clear that the results of Statement 3 and 4 of Theorem \ref{mult} are valid if one of the copulas generates stationary aperiodic Markov chains and another copula of the combination generates stationary ergodic Markov chains.
 
It is also clear that Theorem \ref{mult} holds for a Markov chain generated by the convex combination and any marginal distribution, while Theorem \ref{psi} will require an absolutely continuous marginal distribution. 
 
It is not clear if $\rho^{*}$-mixing can have the property stated in Theorem \ref{mixture} for $\beta$-mixing, $\rho$-mixing, $\psi^{\prime}$-mixing and $\phi$-mixing. We do not have evidence that it is not the case, but we don't yet have the tools to investigate the case.
\end{conclusion}

\begin{acknowledgment}

The author thanks the anonymous referees for valuable comments that have helped transform this work to its present form.
The 2014 Summer Research Grant from the College of Liberal Arts of the University of Mississippi has supported this work.
\end{acknowledgment}

\end{document}